\documentclass[11pt]{article}

\usepackage{amscd,amsmath, amssymb}

\numberwithin{equation}{section}

\def\eqref#1{(\ref{#1})}

\newcommand{\arrow}{{\:\longrightarrow\:}}
\newcommand{\Z}{{\Bbb Z}}
\newcommand{\C}{{\Bbb C}}
\newcommand{\R}{{\Bbb R}}

\def\1{\sqrt{-1}\:}

\newcommand{\cac}{{\cal C}}

\renewcommand{\tilde}{\widetilde}
\renewcommand{\bar}{\overline}
\renewcommand{\phi}{\varphi}
\renewcommand{\epsilon}{\varepsilon}

\renewcommand{\leq}{\leqslant}

\newcommand{\Lie}{\operatorname{Lie}}

\newcommand{\Tot}{\operatorname{Tot}}

\newcounter{Mycounter}[section]
\newcounter{lemma}[section]
\renewcommand{\thelemma}{{Lemma \thesection.\arabic{lemma}}}
\newcommand{\lemma}{%
     \setcounter{lemma}{\value{Mycounter}}
     \refstepcounter{lemma}
     \stepcounter{Mycounter}
     {\noindent \bf \thelemma.\ }}

\newcounter{claim}[section]
\renewcommand{\theclaim}{{Claim \thesection.\arabic{claim}}}
\newcommand{\claim}{%
     \setcounter{claim}{\value{Mycounter}}
     \refstepcounter{claim}
     \stepcounter{Mycounter}
     {\noindent \bf \theclaim.\ }}

\newcounter{sublemma}[section]
\newcounter{corollary}[section]
\newcounter{theorem}[section]
\renewcommand{\thetheorem}{{Theorem \thesection.\arabic{theorem}}}
\newcommand{\theorem}{%
     \setcounter{theorem}{\value{Mycounter}}
     \refstepcounter{theorem}
     \stepcounter{Mycounter}
     {\noindent \bf \thetheorem.\ }}

\newcounter{conjecture}[section]
\renewcommand{\theconjecture}{{Conjecture
\thesection.\arabic{conjecture}}}
\newcommand{\conjecture}{%
     \setcounter{conjecture}{\value{Mycounter}}
     \refstepcounter{conjecture}
     \stepcounter{Mycounter}
     {\noindent \bf \theconjecture.\ }}

\newcounter{proposition}[section]
\renewcommand{\theproposition}
       {{Proposition \thesection.\arabic{proposition}}}
\newcommand{\proposition}{%
     \setcounter{proposition}{\value{Mycounter}}
     \refstepcounter{proposition}
     \stepcounter{Mycounter}
     {\noindent \bf \theproposition.\ }}

\newcounter{definition}[section]
\renewcommand{\thedefinition}
       {{Definition~\thesection.\arabic{definition}}}
\newcommand{\definition}{%
     \setcounter{definition}{\value{Mycounter}}
     \refstepcounter{definition}
     \stepcounter{Mycounter}
     {\noindent \bf \thedefinition.\ }}

\newcounter{example}[section]
\renewcommand{\theexample}{{Example \thesection.\arabic{example}}}
\newcommand{\example}{%
     \setcounter{example}{\value{Mycounter}}
     \refstepcounter{example}
     \stepcounter{Mycounter}
     {\noindent \bf \theexample:\ }}

\newcounter{remark}[section]
\renewcommand{\theremark}{{Remark \thesection.\arabic{remark}}}
\newcommand{\remark}{%
     \setcounter{remark}{\value{Mycounter}}
     \refstepcounter{remark}
     \stepcounter{Mycounter}
     {\noindent \bf \theremark.\ }}

\newcounter{problem}[section]
\newcounter{question}[section]
\makeatletter

\newcommand{\ps@verbit}{%
  \renewcommand{\@oddhead}{%
          \scriptsize
          {Einstein-Weyl structures on Vaisman manifolds}
          \hfil\tiny {L. Ornea and M. Verbitsky, June 13, 2006}}
  \renewcommand{\@evenhead}{\@oddhead}
  \renewcommand{\@oddfoot}{\hfil\thepage\hfil}
  \renewcommand{\@evenfoot}{\@oddfoot}}

\@addtoreset{equation}{section}
\@addtoreset{footnote}{section}
\makeatother

\def\blacksquare{\hbox{\vrule width 5pt height 5pt depth 0pt}}
\def\endproof{\blacksquare}

\begin{document}
\begin{center}
{\LARGE\bf
Einstein-Weyl structures  on complex\\[3mm] manifolds
and conformal version \\[3mm]of Monge-Amp\`ere equation.
\\[3mm]
}

Liviu Ornea and Misha Verbitsky\footnote{Misha Verbitsky is an EPSRC
advanced
fellow supported by  EPSRC grant
GR/R77773/01.


}

\end{center}

{\small
\hspace{0.15\linewidth}
\begin{minipage}[t]{0.7\linewidth}
{\bf Abstract} \\
A Hermitian Einstein-Weyl manifold is a complex manifold
admitting a Ricci-flat K\"ahler covering $\tilde M$,
with the deck transform acting on $\tilde M$ by
homotheties. If compact, it admits a canonical 
Vaisman metric, due to Gauduchon. 
We show that a Hermitian Einstein-Weyl structure
on a compact complex manifold is determined
by its volume form. This result is
a conformal analogue of Calabi's theorem
stating the uniqueness of K\"ahler metrics
with a given volume form in a given K\"ahler
class. We prove that
a solution of a conformal version of complex
Monge-Amp\`ere equation is unique.
We conjecture that a Hermitian Einstein-Weyl
structure on a compact complex manifold
is unique, up to a holomorphic automorphism,
and compare this conjecture to Bando-Mabuchi 
theorem.
\end{minipage}
}

\tableofcontents


\section{Introduction}
\label{_Intro_Section_}

\subsection{Calabi-Yau theorem and Monge-Amp\`ere
equations}

S.-T. Yau \cite{_Yau:Calabi-Yau_} has shown 
that a compact manifold  of K\"ahler type
with vanishing first Chern class
admits a unique K\"ahler-Einstein
metric in a given K\"ahler class.
This result is known as Calabi-Yau theorem;
see \cite{besse} for details and implications
of this extremely important work.
Such a metric is called now {\bf Calabi-Yau metric}.
This theorem was conjectured by
E. Calabi (\cite{_Calabi:KE_}),
who also proved that Calabi-Yau metric
is unique, in a a given K\"ahler class.

The idea of the proof was suggested by Calabi,
who has shown that existence of Calabi-Yau metric
is implied by the following theorem, which is true for all compact
K\"ahler manifolds. 

\hfill

\theorem
(\cite{_Yau:Calabi-Yau_}) Let $M$ be a compact 
K\"ahler manifold, $[\omega]\in H^{1,1}(M)$
a K\"ahler class, and $V\in \Lambda^{n,n}(M)$
a nowhere degenerate volume form.
Then there exists a unique K\"ahler form
$\omega_1\in [\omega]$, such that 
\begin{equation} \label{_form_with_prescri_volume_Equation_}
\omega_1^n=\lambda V,
\end{equation}
where $\lambda$ is a constant.

\endproof

\hfill

Given two K\"ahler forms $\omega$, $\omega_1$, 
in the same K\"ahler class, we can always have
\[
\omega_1= \omega+ d d^c \phi
\]
for some function $\phi$ (this statement
is a consequence of the famous {\em $dd^c$-lemma};
see e.g. \cite{_Griffi_Harri_}).
Then \eqref{_form_with_prescri_volume_Equation_} becomes
\begin{equation} \label{_Monge-Amp_Equation_}
(\omega+d d^c \phi) ^n=\lambda e^f \omega^n.
\end{equation}
Here the function $\phi$ is an unknown, 
the K\"ahler form $\omega$ and the
function $f$ are given, and the 
constant $\lambda$ is expressed through
$f$ and $\omega$ as
\[
\int_M \omega^n = \lambda \int_M e^f \omega^n.
\]

The equation \eqref{_Monge-Amp_Equation_}
is called {\bf the complex Monge-Amp\`ere equation}.
Solutions of \eqref{_form_with_prescri_volume_Equation_}
are unique on a compact manifold, as shown by Calabi
(\cite{_Calabi:KE_}).

\subsection{Monge-Amp\`ere
equations on LCK-manifolds}

In this note, we generalize the Monge-Amp\`ere equation
to locally conformally K\"ahler geometry, and show
that its solutions are unique.

Recall that a locally conformally K\"ahler (LCK) manifold is a 
complex manifold admitting a  K\"ahler covering 
$\tilde M$, with the deck transform acting on $\tilde M$ 
by holomorphic homotheties. If $\tilde M$ is, in addition, Ricci-flat, 
$M$ is called Hermitian Einstein-Weyl, or 
locally conformally K\"ahler Einstein-Weyl.\footnote{Normally,
one defines Hermitian Einstein-Weyl differently, and then
this definition becomes a theorem; see
\ref{_Einstein_Weyl_R_flat_cove_Claim_}.}
The Hermitian Einstein-Weyl geometry is a conformal
analogue of Ka\"hler-Einstein geometry. 

Since the deck transform group acts on $\tilde M$
conformally, the LCK-\-struc\-ture defines a conformal class
of Hermitian metrics on $M$. A metric in this class is
called {\bf an LCK-metric}. In the literature, the
distinction between ``LCK-metrics'' and ``LCK-structures''
is often ignored. Any compact LCK-manifold $M$ is naturally 
equipped with a special Hermitian metric $\omega$ in its conformal
class, called {\bf the Gauduchon metric} (see \ref{_Gaud_metrics_Definition_}).
It is defined uniquely, up to a constant multiplier.
If, in addition, $\omega$ is preserved by a holomorphic
flow, which acts on $\tilde M$ by homotheties,
the LCK-manifold $M$ is called a {\bf Vaisman}
(this notion was first introduced as {\bf generalized Hopf} 
manifold, but the name proved to be inapropriate)\footnote{Traditionally, 
Vaisman manifolds
are defined differently, and this definition 
becomes a theorem; see \ref{def_vai}.}.
P. Gauduchon has proven that all 
compact Hermitian Einstein-Weyl
manifolds are Vaisman (\cite{gau}). 

Given an LCK-manifold $(M, I, \omega)$,
its Hermitian form $\omega$ satisfies 
$d\omega = \omega\wedge \theta$,
where $\theta$ is a closed 1-form,
called {\bf the Lee form of $M$} (see 
Section \ref{_Vaisman_Section_}). Its
cohomology class $[\theta]$ is called {\bf the Lee class
of $M$}. 

The Monge-Amp\`ere equations on LCK-manifolds
are equivalent to finding a Gauduchon
metric with a prescribed volume. 
Uniqueness of such metric (hence,
uniqueness of solutions of 
locally conformally K\"ahler Monge-Amp\`ere)
is due to the following theorem.

\hfill

\theorem\label{main}
Let $(M,J)$ be a compact complex manifold
admitting a Vaisman structure, and
$V\in \Lambda^{n,n}(M)$ a nowhere degenerate,
positive volume form.  Then $M$ admits at most one 
Vaisman structure with the same Lee class, 
such that the volume form of the corresponding
Gauduchon metric is equal to $V$.

\hfill

We give an introduction to LCK-geometry
in Section \ref{_Vaisman_Section_}, and
explain the properties of Einstein-Weyl structures in Section
\ref{_EW_Section_}. We prove
\ref{main} in Section \ref{_Main_proof_Section_}.
The existence result (not proven in this paper)
should follow from the same arguments, added
to those used for the proof of Calabi-Yau theorem.
In the last section, we explain how other results
oh K\"ahler geometry (the Calabi conjecture and
the Bando-Mabuchi theorem) generalize
to conformal setting.

\hfill

\remark
For a K\"ahler manifold, the metric is uniquely determined by
the volume form and the K\"ahler class in cohomology. 
In a conformal setting, the Vaisman metric is defined
uniquely by the volume and the Lee class.
This happens because a relevant cohomology group is $H^2(M, L)$,
where $L$ is the weight bundle of the conformal structure
(see \ref{_weight_bun_Definition_}). It is easy to 
show that all cohomology of the local system $L$
vanish, cf. \cite[Remark 6.4]{or}.

\subsection{Sasaki-Einstein manifolds and Einstein-Weyl geometry}

The compatibility between a complex structure and a Weyl structure
naturally leads to the LCK-condition. This was observed by I.
Vaisman (see also \cite{pps}). Moreover, as shown by P. Gauduchon
(\cite{gau}), a compact Einstein-Weyl
locally conformally K\"ahler manifold is necessarily Vaisman
(see \ref{_Gauduchon_Theorem_}). Then  \ref{main}
is translated into the uniqueness of an
Einstein-Weyl Vaisman metric with a 
prescribed volume form and the Lee class.

The Vaisman manifolds are intimately related to Sasakian
geometry (see \emph{e.g.} \cite{ov1}). Given a Sasakian manifold
$X$, the product $S^1\times X$ has a natural Vaisman structure.
Conversely, any Vaisman manifold admits a canonical
Riemannian submersion to $S^1$, with fibers which
are isometric and equipped with a natural Sasakian structure.

Under this correspondence, the Einstein-Weyl Vaisman
manifolds correspond to Sasaki-Einstein manifolds.
The Sasaki-Einstein manifolds
recently became a focus of much research, due to
a number of new and unexpected examples
constructed by string physicists
(see \cite{_Martelli_Sparks_Yau_},
\cite{_CLPP_}, \cite{_GMSW1_},  \cite{_GMSW2_},
and the references therein).  For a physicist, Sasaki-Einstein
manifolds are interesting because of AdS/CFT correspondence in
string theory. From the mathematical point
of view, these examples are as mysterious as
the Mirror Symmetry conjecture 15 years ago.

The Sasakian manifolds, being transverse
K\"ahler\footnote{This viewpoint was systematically developed in the work of C.P. Boyer, K. Galicki and collaborators. See \emph{e.g.} \cite{_BG:Book_}.}, 
can be studied by the means of algebraic
geometry. One might  hope to obtain and study the Sasaki-Einstein
metrics by the same kind of procedures as used to
study the K\"ahler-Einstein metrics in algebraic
geometry. However, this analogy is not perfect.
In particular, it is possible to show that the
Sasaki-Einstein structures on CR-manifolds are
not unique. 

One may hope to approach the classification of Sasaki-Einstein
structures using the Einstein-Weyl geometry.


\section{Vaisman manifolds and LCK-geometry}
\label{_Vaisman_Section_}


We first review the necessary notions of locally conformally K\"ahler
geometry. See \cite{drag}, \cite{ov1},
\cite{ov2}, \cite{ov3}, \cite{_Verbitsky:LCHK_} for details and examples.

Let $(M,J,g)$ be a complex Hermitian manifold of complex dimension $n$.
Denote by $\omega$ its fundamental two-form $\omega(X,Y)=g(X,JY)$.

\hfill

\definition 
A Hermitian metric $g$ on 
$(M,J)$ is {\bf locally conformally K\"ahler} (LCK for short) if
$$d\omega=\theta\wedge\omega,$$
for a closed 1-form $\theta$.

\hfill

Clearly, for any smooth function $f:\; M \arrow \R^{>0}$, 
$f\omega$ is also an LCK-metric. A conformal class of
LCK-metrics is called {\bf an LCK-structure}.

\hfill

The form $\theta$ is called {\bf the Lee form of the LCK-metric},
 and the dual vector field
$\theta^\sharp$ is called {\bf the Lee field}.

The one-form $\frac 1 2 \theta$ can be interpreted as a (flat) connection
one-form in the bundle of densities of weight
1, usually denoted $L$. This is the real line bundle associated to
the representation
\[
    A\mapsto \mid \det(A)\mid^{\frac{1}{2n}}, \ \ A\in \mathrm{GL}(2n,\R)
\] 

\hfill

\definition\label{_weight_bun_Definition_}
The bundle $L$, equipped with a connection $\nabla_0+\theta$,
is called {\bf the weight bundle of a locally conformally
K\"ahler structure}. One could consider the form $\omega$
as a closed, positive $(1,1)$-form, taking values in $L^2$.

\hfill

\remark
Passing to a covering, we may assume that the flat bundle $L$
is trivial. Then $\omega$ can be considered as a closed,
positive $(1,1)$-form taking values in a trivial vector bundle,
that is, a K\"ahler form. Therefore, any LCK-manifold
admits a covering $\tilde M$ which is K\"ahler. The deck transform
acts on $\tilde M$ by homotheties. This property can be used
as a definition of LCK-structures (see Section \ref{_Intro_Section_}).

\hfill

\definition
A {\bf Vaisman manifold} is an LCK-manifold
equipped with an LCK-metric $g$ whose Lee form is
parallel with respect to the Levi-Civita connection of
$g$. In this case, $g$ is called {\bf a Vaisman metric}.

\hfill

A Vaisman metric on a compact manifold
is unique, up to a constant multiplier.
The proof is due to P. Gauduchon 
(\cite{_Gauduchon_1984_}). 

\hfill

\definition\label{_Gaud_metrics_Definition_}
Let $M$ be an LCK-manifold, $g$ an LCK-metric,
and $\theta$ the corresponding Lee form. The metric
$g$ is called {\bf Gauduchon metric} if $d^*\theta=0$.

\hfill

\theorem\label{_Gauduchon_metric_unique_Theorem_}
(\cite{_Gauduchon_1984_})
Let $M$ be a compact  LCK-manifold. 
Then $M$ admits a Gauduchon metric,
which is unique, up to a constant multiplier.

\endproof

\hfill

\remark
A Vaisman metric is obviously Gauduchon. Indeed,
the Vaisman condition $\nabla\theta=0$
implies the Gauduchon condition $d^*\theta=0$.
Therefore, \ref{_Gauduchon_metric_unique_Theorem_}
implies the uniqueness of Vaisman metrics.

\hfill

\definition
Let  $(\cac, g, \omega)$ be a K\"ahler manifold.
Assume that $\rho$ is a free, proper action
of $\R^{>0}$ on $\cac$, and $g$ and $\omega$ are homogeneous
of weight 2:
\[
\Lie_v\omega= 2\omega, \ \  \Lie_v g= 2g,
\]
where $v$ is the tangent vector field of $\rho$.
The quotient $\cac/\rho$ is called {\bf a Sasakian manifold}.
If $N=\cac/\rho$ is given, $\cac$ is called {\bf the K\"ahler cone
  of $N$}. As a Riemannian manifold, $\cac$ is identified
with the {\bf Riemannian cone} of $(N,g_N)$, 
$\cac(N)=(N\times \R^{>0}, t^2 g_N + dt^2)$.

The Sasakian manifolds are discussed in \cite{_BG:Book_},
in great detail.

\hfill

The following characterization of \emph{compact} Vaisman manifolds
is known (see \cite{ov1}):

 \hfill

\remark\label{def_vai} 
A compact complex manifold $(M,J)$ is
 Vaisman
if it admits a K\"ahler covering $(\tilde
M, J, h)\rightarrow (M,J)$ such that:
\begin{itemize}
\item The monodromy group
$\Gamma\cong \Z$ acts on $\tilde M$ 
by holomorphic homotheties with respect to $h$
(this means that $(M,J)$ is equipped with an LCK-structure).
\item $(\tilde M, J, h)$ is isomorphic to a K\"ahler cone over a
compact Sasakian manifold $S$. Moreover, there exists a Sasakian
automorphism $\phi$ and a positive number $q>1$ such that $\Gamma$
is isomorphic to the cyclic group generated by $(x,t)\mapsto
(\phi(x), tq)$.
\end{itemize}

\hfill

\remark \label{sasa_reeb} 
In these assumptions, denote by $\theta^\sharp$ the vector 
field $t \frac d{dt}$ on 
$\tilde M = (S\times \R^{>0}, g_S t^2 + dt^2)$.
Chose the metric $g= g_S + dt^2$ on $M= \tilde M/\Gamma$.
Clearly, $\theta^\sharp$ descends to a Lee field
on $M$, denoted by the same letter. Then
$J(\theta^\sharp)$ is tangent to the fibers
of the natural projection $\tilde M\arrow \R^{>0}$,
hence belongs to $TS$. This vector field is
called {\bf the Reeb field} of the Sasakian
manifold $S$. Clearly, the orbits of $J(\theta^\sharp)$
on $\tilde M$ are precompact (contained in a compact set).

\hfill

\remark \label{_pote_on_cove_of_Vaisman_Remark_}
It will be important for us to note that the K\"ahler metric $h$ on
the covering $\tilde M =\cac(S)=S\times \R^{>0}$
has a global K\"ahler potential $\phi$, which is
expressed as $\phi(x,t) =t^2$. The
metric $\phi^{-1}\cdot h$ projects on $M$ into the LCK metric
$g$.  Therefore, the Vaisman Hermitian form $\omega$ of $M$ is related
to the K\"ahler form $\tilde \omega$ on $\tilde M$
as follows:
\begin{equation}\label{_Vaisman_via_Kah_and_pot_Equation_}
\omega= \phi^{-1} \tilde \omega
\end{equation}
Moreover, $\phi=\mid\theta\mid^2$, the norm
being taken with respect to the lift of $g$.

\hfill

On a Vaisman manifold, the Lee field $\theta^\sharp$ is Killing,
parallel and holomorphic. One easily proves
that $\mathcal{L}_{\theta^\sharp}\omega=2\omega$.

\hfill

Recall from \cite{tondeur} the notion of transverse geometry:

\hfill

\definition
Consider a manifold endowed with a foliation $\mathcal{F}$
with tangent bundle $F$ and normal bundle $Q$.
A differential, or Riemannian, form
$\alpha$ on $X$ is {\bf basic} (or {\bf transverse}) if  
$X\rfloor\alpha=0$ and $\mathrm{Lie}_X\alpha=0$ for every $X\in F$.
A {\bf
  transverse geometry} of $\mathcal{F}$ is a geometry
defined locally on the leaf space of $\cal F$.
A {\bf K\"ahler transverse structure} on $(M, \cal F)$
is a complex Hermitian structure on $Q$
defined by a  pair $g_{\cal F}, \omega_{\cal F}$ of transverse
forms, in such a way that the induced almost
complex structure defined locally on
the leaf space $M/\cal F$ is 
integrable and K\"ahler.

\hfill

\example\label{_transverse_via_Vaisman_Example_}
Let $(M,J,\omega)$ be a Vaisman manifold,
 $\theta^\sharp$ its Lee field. Consider the
holomorphic foliation ${\cal F}$, generated by
$\theta^\sharp$ and $J\theta^\sharp$.
The form $\omega-\theta\wedge J\theta$ is 
transverse K\"ahler.
Hence the Vaisman manifolds provide examples of transverse K\"ahler
foliations (\cite{_Vaisman_},
\cite{_Tsukada_1994_}). Similarly, a Sasakian manifold has
a transverse K\"ahler geometry associated to the foliation
generated by the Reeb field.

\hfill

A compact complex manifold of Vaisman type
can have many Vaisman structures, still the Lee field is
unique up
to homothety:

\hfill

\proposition\label{tsuk}
If $g_1$, $g_2$ are Vaisman metrics on the same compact manifold
$(M,J)$, then
$\theta_1^{\sharp_{g_1}}=c\theta_2^{\sharp_{g_2}}$,
for some real constant $c$.

\hfill

\noindent{\bf Proof:}
The result was proven by Tsukada in \cite{tsukada}. Here we include
an alternative proof. Recall from \cite{_Verbitsky:LCHK_} that for a
Vaisman structure $(g,J)$, the two-form
\[ \eta:=\omega-\theta\wedge J\theta
\] is exact and positive, with the null-space
generated by  $\langle \theta^\sharp,
J\theta^\sharp\rangle$.
It is the transverse K\"ahler form of $(M, {\cal F})$
(see \ref{_transverse_via_Vaisman_Example_}). 
Let $g_1$, $g_2$ be Vaisman metrics,
$\omega_1$, $\omega_2$ the corresponding Hermitian forms,
$\theta_i$ and $\theta_i^\sharp$ the corresponding
Lee forms and Lee fields.
Consider the $(1,1)$-forms
$\eta_1$, $\eta_2$, defined as above,
\[
\eta_i:=\omega_i-\theta_i\wedge J\theta_i.
\]
Unless their null-spaces coincide,
the sum $\eta_1+\eta_2$ is strictly positive.
Then
\[
  \int_M(\eta_1+\eta_2)^{\dim M}>0.
\]
This is impossible, because $\eta_i$ are exact.
We obtained that the 2-dimensional bundles
generated by $\theta_i^\sharp, J\theta_i^\sharp$ are equal:
\[ \langle \theta_1^\sharp, J\theta^\sharp_1\rangle =
   \langle \theta_2^\sharp, J\theta^\sharp_2\rangle
\]
This implies that $\theta_1^\sharp$,
considered as a vector in $T^{1,0}(M)$, is proportional
to $\theta_2^\sharp$ over $\C$.
\begin{equation}\label{_theta_i_propo_Equation_}
\theta_1^\sharp= a \theta_2^\sharp+ bJ\theta^\sharp_2, \ \ a, b\in \R.
\end{equation}
Since $\theta_i^\sharp$ is holomorphic, the proportionality
coefficient is constant.

To finish the proof of
\ref{tsuk}, it remains to show that this
proportionality coefficient is real.
Here we use \ref{sasa_reeb}: the orbits of $J\theta_1^\sharp$ should be pre-compact.
{}From \eqref{_theta_i_propo_Equation_} we obtain
\[
J \theta_1^\sharp= a J\theta_2^\sharp- b\theta^\sharp_2.
\]
But $aJ\theta_2^\sharp - b \theta^\sharp_2$ acts
on the metric by a homothety, with a coefficient
which is proportional to $e^{-b}$. Therefore,
an orbit of this vector field
is contained in a compact set if and only if $b=0$.

\endproof

\hfill

\remark 
Let $L_\C=L\otimes_\R\C$ be the complexification of the weight
bundle of the Vaisman manifold
$(M,J,g)$. The Lee form then is the connection form of the
standard Hermitian connection in $L_\C$, and one can prove (see
\cite{_Verbitsky:LCHK_}) that its curvature
can be identified with the above form
$\eta=\omega-\theta\wedge J\theta$, 
hence it is exact.


\section{Monge-Amp\`ere equation on LCK-manifolds}
\label{_Main_proof_Section_}


In this Section, we prove \ref{main}. Clearly,
 \ref{main} is implied by the following proposition. 

\hfill

\proposition\label{_omegas_equal_from_det_Proposition_}
Let $(M,J)$ be a compact complex $n$-dimensional manifold admitting two Vaisman
metrics $\omega_1$ and $\omega_2$, such that the
corresponding Lee classes and the volume forms are equal:
\begin{equation}\label{_top_powers_equai_Equation_}
[\theta_1]=[\theta_2], \ \ \  \omega_1^n=\omega_2^n.
\end{equation}
Then $\omega_1=\omega_2$.

\hfill

\noindent{\bf Proof:} We start with the following claim,
which is implied by Tsukada's theorem (\ref{tsuk}).

\hfill

\claim \label{_Lee_fields_equal_Claim_}
In these assumptions, denote the corresponding
Lee fields by $\theta^\sharp_i$, $i=1,2$. Then
\[ \theta^\sharp_1=\theta^\sharp_2.
\]

\hfill

\noindent {\bf Proof:}
Denote by $\tilde \omega_i$ the K\"ahler forms on $\tilde M$
corresponding to $\omega_i$.
By construction, $\Lie_{\theta^\sharp_i}\tilde \omega_i=2\tilde \omega_i$,
where $\Lie$ denotes the Lie derivative.
Therefore,
\begin{equation}\label{_diffe_tilde_ome_Equation_}
\Lie_{\theta^\sharp_i}\tilde \omega_i^n=2n\tilde \omega_i^n.
\end{equation}
Since $[\theta_1]=[\theta_2]$, the forms $\tilde
\omega_1$, $\tilde\omega_2$ have the same automorphy
factors under the deck transform of $\tilde M$,
hence the function $\frac{\tilde\omega_1^n}{\tilde\omega_2^n}$
is invariant under the monodromy of $\tilde M$.
We shall consider $\Psi:=\frac{\tilde\omega_1^n}{\tilde\omega_2^n}$
as a function on $M$. By \ref{tsuk}, $\theta^\sharp_1=c\theta^\sharp_2$.
From \eqref{_diffe_tilde_ome_Equation_}, we obtain
\begin{equation}\label{_Lee_of_frac_Equation_}
\Lie_{\theta^\sharp_1}\Psi= 2n (1-c)\Psi. 
\end{equation}
The function $\Psi$ and the vector field $\theta_1^\sharp$
are defined on a compact manifold $M$; but, unless $c=1$,
\eqref{_Lee_of_frac_Equation_} implies that
the maxumum of $\Psi$ decreases or increases monotonously
under the action of the coresponding flow of diffeomorphisms.
Therefore, \eqref{_Lee_of_frac_Equation_} is possible only
if $c=1$. We proved \ref{_Lee_fields_equal_Claim_}.
\endproof

\hfill

Return to the proof of 
\ref{_omegas_equal_from_det_Proposition_}. Consider a form
\begin{equation}\label{_eta_i_defini_Equation_} 
   \eta_i:=\omega_i-\theta_i\wedge J\theta_i.
\end{equation}
This is a positive, exact $(1,1)$-form on $M$,
which can be interpreted as a curvature of the weight bundle
(see the proof of \ref{tsuk}). 
First of all, we deduce from $\eta_1=\eta_2$ 
the statement of \ref{_omegas_equal_from_det_Proposition_}.

\hfill

\lemma\label{_eta_1=_eta_2_then_omega_equal_Lemma_}
In the assumptions of \ref{_omegas_equal_from_det_Proposition_},
assume that $\eta_1=\eta_2$, where $\eta_i$ are $(1,1)$-forms
defined in \eqref{_eta_i_defini_Equation_}. Then
$\omega_1=\omega_2$.

\hfill

\noindent{\bf Proof.} As follows from \eqref{_eta_i_defini_Equation_},
to prove $\omega_1=\omega_2$ it suffices to show
$\theta_1=\theta_2$.
Let $\tilde M$ be the K\"ahler $\Z$-covering of $M$,
which is a cone over a compact Sasakian manifold,
and $\phi_1, \phi_2$ the corresponding K\"ahler potentials,
obtained as in \ref{_pote_on_cove_of_Vaisman_Remark_}.
It is easy to see that $\theta_i = d \log \phi_i$
and $\eta_i = d^c \theta_i$ (\cite{_Verbitsky:LCHK_}).
Therefore,
\begin{equation}\label{_eta_i_diff_via_d_d^c_Equation_}
\eta_1-\eta_2 = d^c d \log\left(\frac{\phi_1}{\phi_2}\right) 
\end{equation}
The functions $\phi_i$ are automorphic under
the deck transform action on $\tilde M$,
with the same factors of monodromy. Therefore,
their quotient $\frac{\phi_1}{\phi_2}$ is well 
defined on $M$. By
\eqref{_eta_i_diff_via_d_d^c_Equation_},
$0=\eta_1-\eta_2 = d^c d
\log\left(\frac{\phi_1}{\phi_2}\right)$, hence
$\psi:=\log\frac{\phi_1}{\phi_2}$ is pluriharmonic on a compact
complex manifold $M$. Therefore $\psi$ is constant.
This gives $\theta_1-\theta_2 = d\psi=0$.
\ref{_eta_1=_eta_2_then_omega_equal_Lemma_} is proven.
\endproof

\hfill

Return to the proof of
\ref{_omegas_equal_from_det_Proposition_}.
Note that $\eta_i$ are transverse
K\"ahler forms. Since
\[ 
  \det \eta_i = 
  (\theta^\sharp\wedge J\theta^\sharp)\rfloor\det\omega_i,
\]
it follows that
\[ \det\eta_1=\det\eta_2.\]

Let $\rho$ be a transverse form, defined as
$\rho=\sum_{k+l=n-2}\eta_1^k\wedge\eta_2^l$. Then
\begin{equation}\label{_diffe_top_powers_Equation_}
  (\eta_1-\eta_2)\wedge\rho=0.
\end{equation}

As $\eta_i$ are both positive, $\rho$ is
strictly positive, transversal $(n-2, n-2)$-form.
It is well known that on a complex manifold $X$,
any positive $(\dim X-1, \dim X-1)$-form is
an $(\dim X-1)$-st power of a Hermitian form. Therefore,
there exists a transverse form $\alpha$
such that $\rho=\alpha^{n-2}$. Then \eqref{_diffe_top_powers_Equation_} gives
\[
   (\eta_1-\eta_2)\wedge \alpha^{n-2}=0.
\]

{}From \eqref{_eta_i_diff_via_d_d^c_Equation_}, we obtain
$$\eta_1-\eta_2=dd^c\psi,$$
where $\psi:=\log\left(\frac{\phi_1}{\phi_2}\right)$
is a smooth, transversal function on $M$.

We now associate to $\alpha$ a second-order 
differential operator $\mathcal{D}$ acting on transverse $\mathcal{C}^\infty$
functions, which is defined as follows. For any transverse function
$f$,  $dd^c f\wedge \alpha^{n-2}$ is a
transverse top $(n-1, n-1)$ form, and hence there
exists a unique transverse function $g$ such that
$dd^cf\wedge \alpha^{n-2}=g\cdot \alpha^{n-1}$.
We define
\[ \mathcal{D}(f) = g, \quad \text{where} \,\,\,
dd^c f\wedge \alpha^{n-2}=g\cdot \alpha^{n-1}.
\]
In other words,
\[
\mathcal{D}(f) = \frac{dd^cf\wedge \alpha^{n-2}}{\alpha^{n-1}}.
\]
{}From the definition, we have $D(\psi)=0$.
Obviously $\mathcal{D}$ has positive symbol on
the ring of transverse functions, identified
locally with functions on a space of leaves 
of ${\cal F}$.\footnote{In fact, the symbol of $\cal D$ is
  equal to the Riemannian form
  associated with $\alpha$.}
This allows us to
apply the generalized maximum principle:

\hfill

\proposition
(\cite{prot}) Let $\mathcal{D}$
be a second order differential operator on $\R^n$
with positive symbol, satisfying $\mathcal{D}(const.)=0$,
and let $f\in \ker \mathcal{D}$ be a function in its kernel.
Assume that $f$ has a local maximum. Then $f$ is constant.
\endproof

\hfill

Return to the proof of \ref{main}.
Recall that from \eqref{_eta_i_diff_via_d_d^c_Equation_}, we have
\[
\eta_1-\eta_2 = d^c d \log\left(\psi\right), \ \ \psi\in\ker\mathcal{D}.
\]
To show that
$\eta_1=\eta_2$ it is enough to prove that the
kernel of $\mathcal{D}$ contains only constant functions.
As follows from the generalized maximum principle,
a function in $\ker\mathcal{D}$ which has a local maximum
is necessarily constant. Since $M$ is compact,
any continuous function on $M$ must have a maximum. Therefore,
$\psi\in \ker \mathcal{D}$ is constant,
and $\eta_1-\eta_2=d d^c\psi=0$. 
The proof of \ref{_omegas_equal_from_det_Proposition_}  is
finished.
\endproof

\hfill

\remark\label{_Calabi-Yau_for_LCK_Remark_}
In this Section, we proved the uniqueness
of a Vaisman metric satisfying $\omega^n =V$,
for a given non-degenerate, positive volume
form $V$ on a given compact manifold admitting
a Vaisman structure. Following
the argument used to prove the Calabi-Yau theorem
it seems to be possible to prove {\em existence}
of such solutions as well. Indeed, 
$\omega^n =V$ is equivalent to 
\[ 
  \eta^{n-1}=(\theta^\sharp\wedge
  J\theta^\sharp)\rfloor V, 
\]
and this is a transversal Calabi-Yau
equation, which is implied by the same
arguments as used for the
usual Calabi-Yau theorem.

\hfill

\remark
Assuming that Calabi-Yau-type result
(as in \ref{_Calabi-Yau_for_LCK_Remark_}) is true
in LCK geometry, we obtain a very simple 
description of the moduli of Vaisman metrics.
Indeed, the Vaisman metrics in this case
are in one-to-one correspondence with
non-degenerate, positive volume forms.

\section{Einstein-Weyl LCK manifolds}
\label{_EW_Section_}

Einstein-Weyl structures are defined and studied for their own, see
\emph{e.g.} \cite{calder}. Here we specialize the
definitions to LCK structures.

The Levi-Civita connection $\nabla^g$ of $g$ is not the best tool to
study the conformal properties of an LCK manifold.
Instead, the {\bf Weyl connection} defined by

\[
  \nabla =\nabla^g -\frac 1 2\{\theta\otimes Id+Id\otimes \theta
  + g\otimes \theta^\sharp\}
\]
is torsion-free and satisfies $\nabla g=\theta\otimes g$.

The Ricci tensor of the Weyl connection is 
not symmetric. Hence, to obtain the analogue of the Einstein
condition one gives:

\hfill

\definition 
An LCK-manifold is {\bf Einstein-Weyl}
if the symmetric part of the Ricci tensor of the Weyl connection
is proportional to the metric. An Einstein-Weyl LCK-manifold
is also called {\bf Hermitian Einstein-Weyl}. 

\hfill

Let $\nabla$ be a Weyl connection on an LCK-manifold.
One can see that $\nabla$ is the covariant derivative associated to
the connection one-form $\theta$ in
the weight bundle $L$. Since $\theta$ is closed,
we can take a covering $\tilde M$ of $M$, 
with $\theta=df$, for some function $f$ on $\tilde M$.
The Weyl connection becomes the Levi-Civita connection
for the metric $e^{-f}g$ on $\tilde M$. 
Since $\nabla(e^{-f}g)= \nabla(J)=0$,
$e^{-f}g$ is a K\"ahler metric. This
way one obtains a K\"ahler covering of
an LCK-manifold, starting from a Weyl connection.
The converse construction is also clear:
The Levi-Civita connection on a K\"ahler covering $\tilde M$ 
of  an LCK-manifold $M$ is independent from  homotheties,
hence descends to $M$, and satisfies the
conditions for Weyl connection.

\hfill

This gives the following claim.

\hfill

\claim\label{_Einstein_Weyl_R_flat_cove_Claim_}
Let $\nabla$ be a 
 Weyl connection on a complex Hermitian manifold.
Then $\nabla$ satisfies the Einstein-Weyl condition if and only if 
$\nabla$ is Ricci-flat on the K\"ahler covering of $M$.

\endproof

\hfill

\remark
\ref{_Einstein_Weyl_R_flat_cove_Claim_}
also follows from \ref{canon} (below).
Indeed, a trivialization of the weight bundle
$L_\C$ induces a trivialization of canonical
class $K= L^{-n}_\C$.

\hfill

{}From a deep result of Gauduchon in \cite{gau}, it follows that:

\hfill

\theorem \label{_Gauduchon_Theorem_}
Let $(M,J,g)$ be a compact Einstein-Weyl LCK manifold. Then
the Ricci tensor of the Weyl
connection vanishes identically and the Lee form is parallel. In
particular, $(M,J,g)$ is Vaisman.

\endproof

\hfill

{}From \ref{_Gauduchon_Theorem_}, we obtain that all K\"ahler
coverings of an Einstein-Weyl LCK-manifold are Ricci-flat.
This property can be used as a definition of 
Einstein-Weyl LCK-manifolds. 

\hfill

The locally conformally K\"ahler Einstein-Weyl structures
can be expres\-sed in terms of  the complexified weight bundle.

\hfill

\proposition\label{canon} 
(\cite[Proposition 5.6]{_Verbitsky:LCHK_})
Let $M$ be an Einstein-Weyl LCK-manifold, $K$
its canonical class, $L_\C$ its weight bundle.
Consider $K$, $L_\C$ as Hermitian holomorphic
bundles, with the metrics induced from $M$. Then
$L_\C^n\cong K^{-1}$.

\endproof


\section{Bando-Mabuchi theorem in conformal setting}


The covering of an Einstein-Weyl LCK-manifold
is Ricci-flat. Exploiting the analogy with the
Calabi-Yau theorem, one could hope to infer the uniqueness,
or even existence, of Einstein-Weyl structure on 
a complex manifold admitting a Vaisman structure. 
Unfortunately, this analogy does not work that well.
A group of holomorphic automorphisms of the simplest
Hermitian Einstein-Weyl manifold, a Hopf
surface $H = (\C^2\backslash 0)/ \Z$,
does not act on $H$ by conformal automorphisms.
A closer look at the Monge-Amp\`ere equation 
controlling the Einstein-Weyl condition explains this
problem immediately.

Let  $(M,J,g)$ be an Einstein-Weyl
LCK-manifold equipped with a Vaisman metric, and $\tilde M$
its K\"ahler covering, which trivializes $L$. From \ref{canon},
it is clear that $\tilde M$ has trivial canonical class.
Let $\Omega$ be a section of canonical class of $\tilde M$ which
is equivariant under the monodromy action. Such a section is unique
up to a constant. Indeed, if $\Omega_1$, $\Omega_2$
are two equivariant sections of canonical class,
the quotient $\frac{\Omega_1}{\Omega_2}$
is a holomorphic function on $\tilde M$
which is invariant under monodromy, hence
descends to a global holomorphic function on $M$.
Therefore $\frac{\Omega_1}{\Omega_2}=const$. 
Rescaling $\Omega$  such that $|\Omega|=1$, we obtain
\[
\Omega\wedge\bar\Omega= \frac{1}{n!\,2^n} \tilde \omega^n,
\]
where $n=\dim_\C M$ and $\tilde \omega$ is the
K\"ahler form on $\tilde M$. In particular, given two Einstein-Weyl
structures with the K\"ahler forms
$\tilde\omega_1$ and $\tilde\omega_2$, we always have
$\tilde\omega_1^n = \lambda \tilde\omega_2^n$, where $\lambda$
is a positive constant. After rescaling, we may also
assume that
\begin{equation}\label{_tilde_omega_i^top_independent_Equation_}
\det \tilde\omega_1 = \det \tilde\omega_2,
\end{equation}
where $\det\tilde\omega_i=\tilde \omega_i^n$, $n=\dim_\C M$.
Comparing this equation with \eqref{_Vaisman_via_Kah_and_pot_Equation_},
we find that  \eqref{_tilde_omega_i^top_independent_Equation_}
is equivalent to 
\begin{equation}\label{_omega_i_Hermitian^top_independent_Equation_}
\frac{\det \omega_1}{\det \omega_2}= \frac{\phi_2}{\phi_1},
\end{equation}
where $\phi_i$ denotes the K\"ahler potential of the K\"ahler
metrics $\omega_i$. Writing
$\psi_i := \log\phi_i$, $\psi:= \psi_1-\psi_2$,
$\eta_i := d d^c \psi_i$, $\eta :=\eta_1-\eta_2$, and using 
\begin{equation}
\frac{\det \omega_1}{\det \omega_2}= \frac{\det \eta_1}{\det \eta_2}
\end{equation}
(see the proof of \ref{_omegas_equal_from_det_Proposition_}),
we find that 
\eqref{_omega_i_Hermitian^top_independent_Equation_}
is equivalent to the transversal version 
of Aubin-Calabi-Yau equation:
\begin{equation}
\log\left( \frac{\det(\eta- dd^c \psi)}{\det \eta}\right) 
= \epsilon\psi+const,
\end{equation}
with $\epsilon=1$. This equation is well-known in the theory
of K\"ahler-Einstein manifolds. It is easy to see that
it has a unique solution $\psi=0$ for $\epsilon\leq 0$. 
When $\epsilon=1$, this becomes much more difficult.
Bando and Mabuchi (\cite{_Bando_Mabuchi_}) 
studied this equation in order to prove that a 
 K\"ahler-Einstein metric on a Fano manifold
is unique, up to a constant multiplier and a 
holomorphic automorphism. We hope that their
argument will carry over in Einstein-Weyl geometry,
proving the following conjecture.

\hfill

\conjecture\label{_M-B_for_LCK_Conjecture_}
Let $(M,J)$ be a compact complex manifold.
Then it admits at most one
Einstein-Weyl structure, up to a
holomorphic automorphism. 

\hfill

This conjecture cannot be very easy, because it implies
Bando-Mabuchi theorem. Indeed, assume that 
$M= \Tot(K\backslash 0)/\langle q\rangle$,
where $K\backslash 0$ is the space of 
non-zero vectors in the canonical line bundle
of a Fano manifold $X$, and $\langle q\rangle$
a holomorphic $\Z$-action generated by $v \arrow qv$,
where $q$ is a fixed complex number, $|q|>1$.
The manifold $M$ clearly admits a Vaisman
structure; it is Einstein-Weyl if and only if $X$
is K\"ahler-Einstein, and the K\"ahler-Einstein
metrics on $X$ correspond uniquely to
the Einstein-Weyl structures on $M$.
Therefore, \ref{_M-B_for_LCK_Conjecture_}
for such $M$ is equivalent to the
Bando-Mabuchi theorem for $X$.

\hfill


\hfill

\noindent{\bf Acknowledgements:}
This note originated from discussions at the 
Max Planck Institute in Bonn
whose support is gratefully acknowledged. The second author
is grateful to Semyon Alesker, who explained to him the
generalized maximum principle.

\hfill

{\small

\noindent {\sc Liviu Ornea\\
University of Bucharest, Faculty of Mathematics, \\14
Academiei str., 70109 Bucharest, Romania.}\\
\tt Liviu.Ornea@imar.ro, \ \ lornea@gta.math.unibuc.ro

\hfill

\noindent {\sc Misha Verbitsky\\
University of Glasgow, Department of Mathematics, \\15
  University Gardens, Glasgow, Scotland.}\\
{\sc  Institute of Theoretical and
Experimental Physics \\
B. Cheremushkinskaya, 25, Moscow, 117259, Russia }\\
\tt verbit@maths.gla.ac.uk, \ \  verbit@mccme.ru
}

\end{document}